\documentclass{article}
\usepackage{amssymb, amsfonts, amsmath, amsthm}
\usepackage{enumerate}
\usepackage{caption}
\usepackage[colorlinks=true,
citecolor=black,
linkcolor=black,
anchorcolor=black,
filecolor=black,
menucolor=black,
urlcolor=black
]{hyperref}
\usepackage[normalem]{ulem}
\usepackage{cancel}
\usepackage{array}
\usepackage{longtable}

\usepackage{epsfig,lpic,wrapfig}
\usepackage{bm}
\usepackage{multicol}

\usepackage[utf8]{inputenc}
\usepackage[T1,T2A]{fontenc}
\usepackage[russian,
german,
english]{babel}

\usepackage[
sorting=debug,
bibencoding=auto,
backend=biber,
maxnames=10,
autolang=other,
doi=false,
url=false,
style=numeric-comp,
isbn=false]{biblatex}
\usepackage{csquotes}
\AtEveryBibitem{\clearlist{language}}
\AtEveryBibitem{\clearlist{note}}
\renewbibmacro{in:}{}

\DeclareFieldFormat{postnote}{#1}
\DeclareFieldFormat{multipostnote}{#1}
\DeclareFieldFormat{pages}{#1}

\def\thetitle{Graph comparison meets Alexandrov}
\def\theauthors{Nina Lebedeva and Anton Petrunin}
\hypersetup{pdftitle={\thetitle},
pdfauthor={\theauthors}}

\newcommand{\Addresses}{{\bigskip\footnotesize

\noindent Nina Lebedeva,
\par\nopagebreak
 \textsc{St. Petersburg State University, 7/9 Universitetskaya nab., St. Petersburg, 199034, Russia}
\par
\nopagebreak
 \textsc{St. Petersburg Department of V. A. Steklov Institute of Mathematics of the Russian Academy of Sciences, 27 Fontanka nab., St. Petersburg, 191023, Russia}
  \par\nopagebreak
  \textit{Email}: \texttt{lebed@pdmi.ras.ru}

\medskip

\noindent   Anton Petrunin, 
\par\nopagebreak
 \textsc{Math. Dept. PSU, University Park, PA 16802, USA.}
  \par\nopagebreak
  \textit{Email}: \texttt{petrunin@math.psu.edu}
  
}}

\def\parit#1{\medskip\noindent{\it #1}}
\def\parbf#1{\medskip\noindent{\bf #1}}
\def\qeds{\qed\par\medskip}
\def\qedsf{\vskip-6mm\qeds}

\newcounter{thm}[section]

\def\claim#1{\par\medskip\noindent\refstepcounter{thm}\hbox{\bf\boldmath #1.}
\it\ 
}
\def\endclaim{
\par\medskip}
\newenvironment{thm}{\claim}{\endclaim}

\begin{document}

\title{\thetitle}
\author{\theauthors}

\date{}
\maketitle
\begin{abstract}
Graph comparison is a certain type of condition on metric space encoded by a finite graph.
We show that any nontrivial graph comparison implies one of two Alexandrov's comparisons.
The proof gives a complete description of graphs with trivial graph comparisons.
\end{abstract}

\parbf{Preface.}
The notion of graph comparison was introduced in \cite{lebedeva-petrunin-zolotov}.
It was studied further in \cite{toyoda,toyoda2019,lebedeva-petrunin-CBB,lebedeva,lebedeva-petrunin,lebedeva-petrunin-octahedron}.
Let us mention some of the results.
\begin{itemize}
\item Graph comparisons for the tripod and four-cycle capture nonnegative and nonpositive curvature in the sense of Alexandrov; see below.
\end{itemize}
\begin{itemize}
\item Graph comparison for star graphs provides a stronger version of the so-called \emph{Lang--Schroeder--Sturm inequality} \cite{lang-schroeder, sturm, lebedeva-petrunin-CBB}.
\end{itemize}
\begin{itemize}
\item The all-tree comparison gives a metric description of target spaces of submetries from subsets of Hilbert space \cite{lebedeva-petrunin-zolotov}.
\end{itemize}

\noindent
\begin{minipage}
{.80\textwidth}
\begin{itemize}
\item The comparison for the tree on the diagram has tight relation with the \emph{transport continuity property} and the so-called \emph{Ma--Trudinger--Wang condition}~\cite{lebedeva-petrunin-zolotov,ma-trudinger-wang}.
\end{itemize}
\end{minipage}
\hfill
\begin{minipage}{.17\textwidth}
\centering
\vskip-1mm
\includegraphics{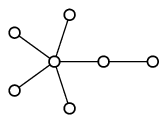}
\end{minipage}
\begin{itemize}
\item Octahedron comparison holds in products of trees~\cite{lebedeva-petrunin-octahedron}.
\end{itemize}

We will show that any nontrivial graph comparison implies one of two Alexandrov's comparisons.

\parbf{Introduction.}
Let us start with the definition.
Suppose $\Gamma$ is a graph with vertices $v_1,\dots,v_n$.
We write $v_i\sim v_j$ (or $v_i\nsim v_j$) if $v_i$ is adjacent (respectively nonadjacent) to $v_j$.

A metric space $X$ meets the \emph{$\Gamma$-comparison} if for any $n$ points in $X$ labeled by vertices of $\Gamma$ there is a model configuration $\tilde v_1,\dots,\tilde v_n$ in the Hilbert space $\mathbb{H}$ such that 
\begin{align*}
v_i\sim v_j\quad&\Longrightarrow\quad|\tilde v_i-\tilde v_j|_{\mathbb{H}}\leqslant | v_i-v_j|_{X},
\\
v_i\nsim v_j\quad&\Longrightarrow\quad|\tilde v_i-\tilde v_j|_{\mathbb{H}}\geqslant | v_i-v_j|_{X};
\end{align*}
here $|\ -\ |_X$ denotes the metric in~$X$.
(Note that $v_i$ refers to a vertex in $\Gamma$ and to the corresponding point in $X$.)

\begin{wrapfigure}{r}{35 mm}
\vskip-0mm
\centering
\includegraphics{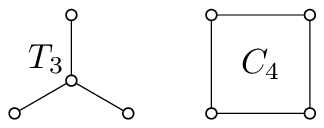}
\end{wrapfigure}

Denote by $T_3$ and $C_4$ the tripod and four-cycle shown on the diagram.
The $C_4$-comparison is equivalent to nonnegative curvature,
and $T_3$-comparison is equivalent to the nonpositive curvature in the sense of Alexandrov \cite{lebedeva-petrunin-zolotov}.
These definitions are usually applied to length spaces, but they can be applied to general metric spaces;
the latter convention is used in \cite{alexander2019alexandrov}.

\begin{thm}{Theorem}
Let $\Gamma$ be an arbitrary finite graph.
Then either $\Gamma$-comparison holds in any metric space,
or it implies $C_4$- or $T_3$-comparison.
\end{thm}

The next statement is a corollary from the proof of the theorem;
it describes all graphs with trivial comparison.

\begin{thm}{Corollary} Let $\Gamma$ be a finite connected graph.
Suppose that $\Gamma$-comparison is trivial;
that is, it holds in any metric space.
Then $\Gamma$ can be constructed from a path $P_{\ell}$ of length $\ell\geqslant 0$ and two complete graphs $K_{m_1}$, $K_{m_2}$ by attaching $k_1$ vertices of $K_{m_1}$ to the left end of $P_{\ell}$ and $k_2$ vertices of $K_{m_2}$ to the right end of~$P_{\ell}$.
\end{thm}

The graph $\Gamma$ in the corollary is described by five integers $(m_1$, $k_1$, $\ell$, $k_2$, $m_2)$ such that $\ell \geqslant 0$, $m_i\geqslant k_i\geqslant 0$, and $k_i>0$ if $m_i>0$ for each $i$.
Examples of such graphs and their 5-arrays are shown below.
\begin{figure}[ht!]
\centering
\includegraphics{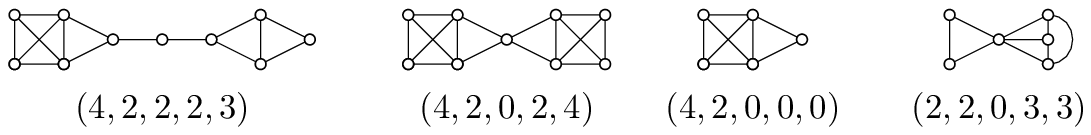}
\end{figure}

\bigskip

\parbf{Proof of the theorem.}
Suppose $\Gamma$ has connected components $\Gamma_1,\dots,\Gamma_k$.
Observe that $\Gamma$-comparison holds in a metric space $X$ if and only if so does every $\Gamma_i$-comparison.
Therefore we can assume that $\Gamma$ is connected.

Suppose $\Gamma$ is a graph with vertices $v_1,\dots,v_n$ as before.
Remove two vertices, say $v_1$ and $v_2$, from $\Gamma$, and add a new vertex $w$ such that for any other vertex $u$ we have
\begin{itemize}
 \item if $u\sim v_1$ and $u\sim v_2$, then $u\sim w$;
 \item if $u\nsim v_1$ and $u\nsim v_2$, then $u\nsim w$;
 \item in the remaining cases, we may choose $u\sim w$ or $u\nsim w$.
\end{itemize}
Denote the so-obtained graph by $\Gamma'$.

Applying the definition of $\Gamma$-comparison assuming that $v_1=v_2$ in $X$, we get the following.

\begin{thm}{Claim}
If $\Gamma$-comparison holds in a metric space $X$, then so does $\Gamma'$-comparison.

\end{thm}

The described construction of $\Gamma'$ from $\Gamma$ will be called \emph{vertex fusion}.
If a graph $\Delta$ can be obtained from $\Gamma$ applying vertex fusion several times, then we will write $\Delta\prec \Gamma$.

From above we get the following two observations:
\begin{itemize}
 \item If $\Delta$ is an induced subgraph of a connected finite graph $\Gamma$, then $\Delta\prec \Gamma$.
 \item If $\Delta\prec \Gamma$, then $\Gamma$-comparison implies $\Delta$-comparison.
\end{itemize}
Hence we get the following reformulation of the theorem.

\begin{thm}{Reformulation}
For any finite connected graph $\Gamma$,
\begin{enumerate}[(a)]
\item $\Gamma$-comparison is trivial,  or
\item $C_4\prec \Gamma$, or
\item $T_3\prec \Gamma$.
\end{enumerate}
\end{thm}

A connected graph will be called \emph{multipath} if it has an integer function $\ell$ on its vertex set such that 
$$v\sim w
\quad\Longleftrightarrow\quad
|\ell(v)-\ell(w)|\leqslant 1.$$
The value $\ell(w)$ will be called the \emph{level} of the vertex $w$.
Multipath is completely described by a sequence of integers that give the number of vertexes on each level.
\begin{figure}[ht!]
\centering
\medskip
\includegraphics{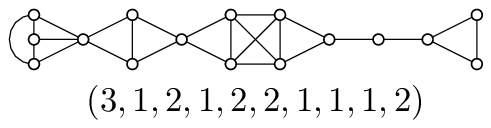}
\medskip
\end{figure}
An example of a multipath with its sequence is shown on the diagram. 

\begin{thm}{Lemma}
Let $\Gamma$ be a connected finite graph such that $C_4\nprec\Gamma$ and  $T_3\nprec\Gamma$.
Then $\Gamma$ is a multipath.
\end{thm}

\parit{Proof.}
We will denote by $|\ -\ |_\Gamma$ the path metric on the vertex set of $\Gamma$;
it is the number of edges in a shortest path connecting two vertices.
Let us show that 
\[|u-w|_\Gamma\geqslant |u-v|_\Gamma\geqslant|v-w|_\Gamma\geqslant 2
\quad\Longrightarrow\quad |u-w|_\Gamma=|u-v|_\Gamma+|v-w|_\Gamma
\leqno({*})\]
for any three vertices $u$, $v$, and $w$ in $\Gamma$.

Suppose $({*})$ does not hold.
Let $\Delta$ be the subgraph of $\Gamma$ induced by three shortest paths between each pair in the triple $u$, $v$, $w$.
\begin{figure}[ht!]
\centering
\includegraphics{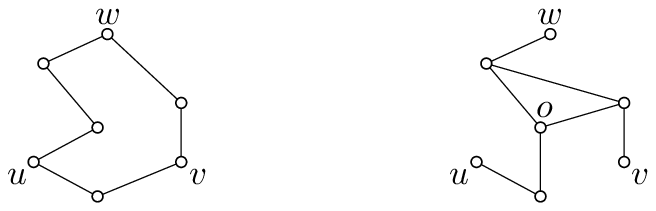}
\end{figure}
Note that $\Delta$ is either a cycle
or it has three paths from a vertex, say $o$, to each of $u$, $v$, and $w$ such that each of these paths does not visit the remaining vertices in the triple.
In these cases, we have $C_4\prec\Delta$ or $T_3\prec\Delta$ respectively.
By the observation above $\Delta\prec\Gamma$; $({*})$ is proved.

Denote by $d$ the diameter of $\Gamma$.
We can assume that $d\geqslant 2$;
if $d=1$, then $\Gamma$ is a complete graph; in particular, it is a multipath.

Choose vertices $p$ and $q$ such that $|p-q|_\Gamma=d$.
Let us show that $\Gamma$ is a multipath with the following level function
\[\ell(w)=
\begin{cases}
\hfil |p-w|_\Gamma&\text{if}\quad |p-w|_\Gamma\geqslant 2,
\\
d-|q-w|_\Gamma&\text{if}\quad |q-w|_\Gamma\geqslant 2,
\\
\hfil 1&\hfil\text{otherwise.}
\end{cases}
\]
By $({*})$, $|p-w|_\Gamma+|q-w|_\Gamma=d$ if $|p-w|_\Gamma\geqslant 2$ and $|q-w|_\Gamma\geqslant 2$;
therefore $\ell$ is well defined.

If $d\geqslant 4$, then the statement follows from $({*})$.
Two cases remain $d=2$ and $d=3$.

\begin{wrapfigure}{r}{17 mm}
\vskip-6mm
\centering
\includegraphics{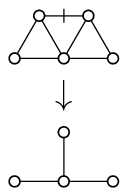}
\vskip-2mm
\end{wrapfigure}

Observe that the fan graph on the diagram cannot appear as an induced subgraph of $\Gamma$.
If this is the case, then applying the vertex fusion to the ends of the marked edge, we could get a tripod --- a contradiction.

\parit{Case $d=2$.}
By $({*})$, $\ell^{-1}(0)$ and $\ell^{-1}(2)$ are cliques.
Observe that $\ell(v)=1$ if and only if $p\sim v$ and $q\sim v$.
Note that $\ell^{-1}(1)$ is a clique;
indeed, if $u\nsim v$ for some $u,v\in\ell^{-1}(1)$,
then the subgraph induced by $\{p,q,u,v\}$ is a four-cycle --- a contradiction.

It remains to show that $u\sim v$, $v\sim w$, and $u\nsim w$ if $\ell(u)=0$, $\ell(v)=1$, and $\ell(w)=2$.

\begin{wrapfigure}{r}{17 mm}
\vskip-2mm
\centering
\includegraphics{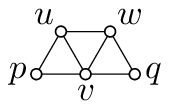}
\vskip4mm
\includegraphics{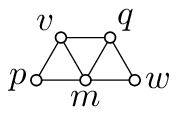}
\end{wrapfigure}

Suppose $u\sim w$. 
Note that $v\sim u$ and $v\sim w$; otherwise, $\Gamma$ contains an induced 4- or 5-cycle with vertices $p, u, v, w, q$.
Therefore the induced subgraph for $\{p,u,v,w,q\}$ is isomorphic to the fan --- a contradiction.

Now, suppose $v\nsim w$; note that $w\ne q$.
Denote by $m$ a midvertex of $w$ and $p$.
From above, $\ell(m)=1$; in particular, $v\sim m\sim q$.
And again, the induced subgraph for $\{p,v,m,w,q\}$ is isomorphic to the fan --- a contradiction.
The same way one shows that $u\sim v$.

\parit{Case $d=3$.}
We need to show that $u\sim v$ if $\ell(u)=2$ and $\ell(v)=3$;
the rest follows from $({*})$.

\begin{wrapfigure}{r}{20 mm}
\vskip-8mm
\centering
\includegraphics{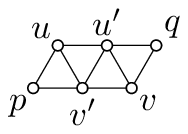}
\vskip-2mm
\end{wrapfigure}

Suppose the contrary;
let $u'$ ($v'$) be a midvertex of $u$ and $q$ (respectively, $v$ and $p$).
Observe that $\ell(u')=3$, $\ell(v')=2$, and $u'\sim v'$ (otherwise the subgraph induced by $\{u,v,u',v'\}$ is a four-cycle).
It follows that the subgraph induced by $\{p,q,u,v,u',v'\}$ is shown on the diagram.
Note that it contains an induced fan --- a contradiction. 
\qeds

\begin{thm}{Proposition}
Let $\Gamma$ be a multipath with sequence $(k_0,\dots, k_m)$.
Suppose $C_4\nprec\Gamma$ and  $T_3\nprec\Gamma$.
Then 
\begin{enumerate}[(a)]
 \item\label{lem:multipath:5} If $m\geqslant 4$, then $k_2=\dots=k_{m-2}=1$.
 \item\label{lem:multipath:4} If $m= 3$, then $k_1=1$ or $k_2=1$.
 \item\label{lem:multipath:3} If $m= 2$, then $k_0=1$, $k_1=1$, or $k_2=1$.
\end{enumerate}

\end{thm}

\parit{Proof.} Assuming the contrary in each case we get

\parit{(\ref{lem:multipath:5})}
if $m\geqslant 4$, then multipath $(1,1,2,1,1)$ is an induced subgraph of~$\Gamma$,

\parit{(\ref{lem:multipath:4})}
if $m=3$, then multipath $(1,2,2,1)$ is an induced subgraph of~$\Gamma$,

\parit{(\ref{lem:multipath:3})}
if $m=2$, then multipath $(2,2,2)$ is an induced subgraph of~$\Gamma$.

In each case, we arrive at a contradiction by applying vertex fusion to the ends of the marked edges as shown on the diagram.

\begin{figure}[h!]
\centering
\includegraphics{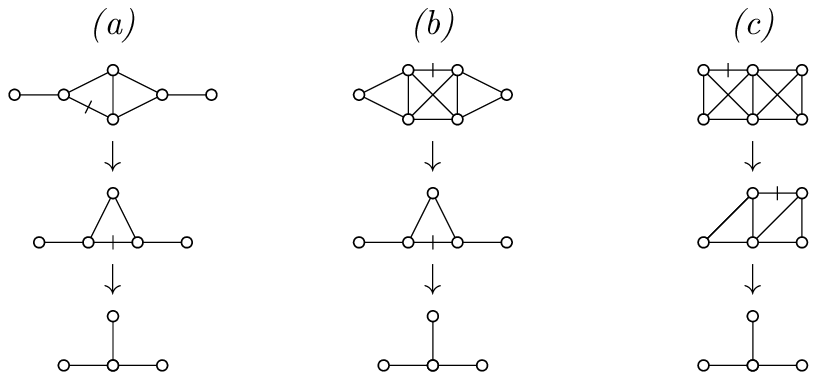}
\end{figure}
\qedsf

It remains to show that $\Gamma$-comparison is trivial for every multipath $\Gamma$ described in the proposition.
This is done by prescribing the coordinates for the needed model configuration on the real line.

Each edge of $\Gamma$ comes with weight --- the distance between the endpoints in~$X$.
Define the distance $\|v-w\|_\Gamma$ as the minimal total weight of paths connecting $v$ to $w$ in $\Gamma$.
Note that 
\[\|v-w\|_\Gamma\geqslant |v-w|_X\]
for any $v$ and $w$.

If $m\leqslant 1$, then $\Gamma$ is a complete graph.
In this case, $\Gamma$-comparison is trivial.
It remains to consider cases $m\geqslant2$.

Let us choose a special vertex $w$ that is unique on its level and not too far from the middle of $\Gamma$.
Namely, if $m\geqslant 4$, then choose $w$ on the second level;
by the proposition, it is unique on its level.
If $m=3$, then by the proposition we can assume that $k_2=1$; in this case choose $w$ on the second level.
Finally, if $m=2$, let $w$ be any vertex that is unique on its level; it exists by the proposition. 

For every vertex $v_i$, let
\[\tilde v_i=\pm \|w-v_i\|_\Gamma,\]
where the sign is plus if $v_i$ has a higher level than $w$ and minus otherwise.
By the triangle inequality, the obtained configuration $\tilde v_1,\dots,\tilde v_n\in\mathbb{R}$ meets the condition of $\Gamma$-comparison.
\qeds

\parbf{Remarks.}
The statements in the preface indicate that for a carefully chosen graph (or a family of graphs) its graph comparison is responsible for meaningful geometric properties of metric spaces.
Let us state two more observations about graph comparison.

Graph comparisons for all complete bipartite graphs imply the so-called \emph{pure inequalities of negative type} \cite[6.1.1]{deza-laurent}.
By Schoenberg's criterion, these inequalities are sufficient for the existence of isometric embedding into a Hilbert space \cite[6.2.1]{deza-laurent}.
In particular, the comparisons for all graphs imply that the metric space is isometric to a subset of a Hilbert space.
The latter statement can be also proved directly the same way as Proposition 1.9 in \cite{toyoda}.

\begin{wrapfigure}{r}{17 mm}
\vskip-6mm
\centering
\includegraphics{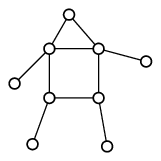}
\vskip-2mm
\end{wrapfigure}

The last observation works for arbitrary metrics.
For length metrics, most graph comparisons imply that the space is isometric to a subset of a Hilbert space.
Indeed if $C_4\prec\Gamma$ and $T_3\prec\Gamma$ (as for the graph on the diagram), then any complete length space that meets $\Gamma$-comparison has vanishing curvature in the sense of Alexandrov; in particular, it is isometric to a convex closed set in a Hilbert space.

\parbf{Acknowledgments.}
We want to thank Alexander Lytchak for help. 

The first author was partially supported by the Russian Foundation for Basic Research, grant 20-01-00070; the second author was partially supported by the National Science Foundation, grant DMS-2005279
and the Ministry of Education and Science of the Russian Federation, grant 075-15-2022-289.

{\sloppy
\printbibliography[heading=bibintoc]
\fussy
}

\Addresses
\end{document}


\title{\thetitle}
\author{\theauthors}

\date{}
\maketitle
\begin{abstract}
Граф-сравнение это условие определённого типа на метрические пространства закодированное конечным графом.
Мы доказываем, что любое нетривиальное граф-сравнение влечёт одно из двух условий Александрова.
Доказательство также даёт классификацию графов с тривиальным сравнением.
\end{abstract}

\parbf{Предисловие.}
Граф-сравнения были введены в \cite{lebedeva-petrunin-zolotov} и изучалось далее в \cite{toyoda,toyoda2019,lebedeva-petrunin-CBB,lebedeva,lebedeva-petrunin,lebedeva-petrunin-octahedron}.
Перечислим некоторые из результатов.
\begin{itemize}
\item Граф-сравнения для треноги и четыре-цикла описывают условия Александрова (см. ниже).
\end{itemize}
\begin{itemize}
\item Граф-сравнения для звёзд дают более сильный вариант так называемых \emph{неравенств Лэнга --- Шрёдера --- Штурма} \cite{lang-schroeder, sturm, lebedeva-petrunin-CBB}.
\end{itemize}
\begin{itemize}
\item Граф-сравнения для всех деревьев даёт описание пространств-образов субметрий из подмножеств гильбертова пространства \cite{lebedeva-petrunin-zolotov}.
\end{itemize}

\noindent
\begin{minipage}
{.80\textwidth}
\begin{itemize}
\item Граф-сравнение для дерева на рисунке, тесно связанно с \emph{непрерывностью решения транспортной задачи} и с так называемым \emph{условием Ма --- Трудингера --- Ванга}~\cite{lebedeva-petrunin-zolotov,ma-trudinger-wang}.
\end{itemize}
\end{minipage}
\hfill
\begin{minipage}{.17\textwidth}
\centering
\vskip-1mm
\includegraphics{mppics/pic-5}
\end{minipage}
\begin{itemize}
\item Октоэдральное сравнение выполняется в произведении деревьев~\cite{lebedeva-petrunin-octahedron}.
\end{itemize}

Мы покажем, что любое нетривиальное граф-сравнение влечёт одно из двух условий Александрова.

\parbf{Введение.}
Начнём с определения.
Пусть $\Gamma$ --- граф с вершинами $v_1,\dots,v_n$.
Будем писать $v_i\sim v_j$ (или $v_i\nsim v_j$) если вершина $v_i$ соединена (соответственно не соединена) ребром с $v_j$.

Метрическое пространство $X$ удовлетворяет \emph{$\Gamma$-сравнению} если для любых $n$ точек в $X$ обозначенных $v_1,\dots,v_n$ существует конфигурация сравнения --- набор точек $\tilde v_1,\dots,\tilde v_n$ в гильбертовом пространстве $\mathbb{H}$ такой, что
\begin{align*}
v_i\sim v_j\quad&\Longrightarrow\quad|\tilde v_i-\tilde v_j|_{\mathbb{H}}\leqslant | v_i-v_j|_{X},
\\
v_i\nsim v_j\quad&\Longrightarrow\quad|\tilde v_i-\tilde v_j|_{\mathbb{H}}\geqslant | v_i-v_j|_{X},
\end{align*}
где $|\ -\ |_X$ обозначает расстояние в~$X$.
(Обратите внимание, что здесь $v_i$ относится и к вершине $\Gamma$ и к соответствующей точке $X$.)

\begin{wrapfigure}{r}{35 mm}
\vskip-0mm
\centering
\includegraphics{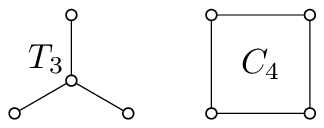}
\end{wrapfigure}

Обозначим через $T_3$ и $C_4$ треногу и четыре-цикл показанные на рисунке.
Их сравнения эквивалентны соответственно неотрицательной и неположительной кривизне в смысле Александрова \cite{lebedeva-petrunin-zolotov}.
Эти определения обычно используются для пространств с внутренней метрикой,
однако их можно применять и к произвольным метрическим пространствам.
Последнее соглашение используется, например, в \cite{alexander2019alexandrov}.

\begin{thm}{Теорема}
Пусть $\Gamma$ --- конечный граф.
Тогда, либо $\Gamma$-сравнение выполняется в любом метрическом пространстве,
любо оно влечёт $C_4$- или $T_3$-сравнение.
\end{thm}

Утверждение ниже является следствием из доказательства теоремы;
оно описывает все графы с тривиальным сравнением.

\begin{thm}{Следствие}
Пусть $\Gamma$ --- конечный связный граф.
Предположим, что $\Gamma$-сравнение тривиально;
то есть, оно выполняется в любом метрическом пространстве.
Тогда граф $\Gamma$ можно построить из пути $P_{\ell}$ длины $\ell\geqslant 0$ и двух полных графов $K_{m_1}$, $K_{m_2}$ подсоединив $k_1$ вершин $K_{m_1}$ к левому концу $P_{\ell}$ и $k_2$ вершин $K_{m_2}$ к правому концу $P_{\ell}$.
\end{thm}

Граф $\Gamma$ в этом следствии определяется пятёркой целых чисел $(m_1$, $k_1$, $\ell$, $k_2$, $m_2)$ таких, что $\ell \geqslant 0$, $m_i\geqslant k_i\geqslant 0$, и $k_i>0$, если $m_i>0$ для всякого $i$. 
Примеры таких графов и их пятёрок показаны ниже.
\begin{figure}[ht!]
\centering
\includegraphics{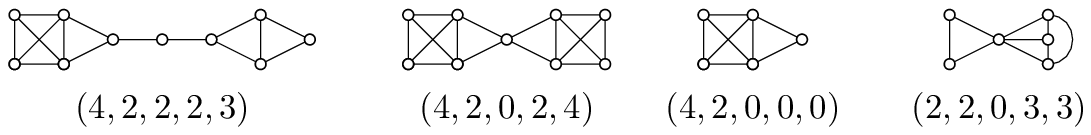}
\end{figure}

\bigskip

\parbf{Доказательство теоремы.}
Предположим, что граф $\Gamma$ распадается на связные компоненты $\Gamma_1,\dots,\Gamma_k$.
Заметим, что $\Gamma$-сравнение эквивалентно выполнению всех $\Gamma_i$-сравнений.
Таким образом можно предположить, что $\Gamma$ связен.

Положим, как и выше, что $v_1,\dots,v_n$ --- вершины $\Gamma$.
Удалим из $\Gamma$ пару вершин, скажем $v_1$ и $v_2$.
Добавим вместо них новую вершину $w$ так, что для любой другой вершины $u$ выполняется следующее
\begin{itemize}
 \item если $u\sim v_1$ и $u\sim v_2$, то $u\sim w$;
 \item если $u\nsim v_1$ и $u\nsim v_2$, то $u\nsim w$;
 \item в остальных случаях можно выбрать $u\sim w$ или $u\nsim w$.
\end{itemize}
Пусть $\Gamma'$ --- полученный граф.

Полагая в определении $\Gamma$-сравнения, что точки $v_1$ и $v_2$ совпадают в $X$,
мы приходим к следующему.

\begin{thm}{Утверждение}
$\Gamma$-сравнение влечёт $\Gamma'$-сравнение.

\end{thm}

Приведённое построение $\Gamma'$ из $\Gamma$ будет называться \emph{слиянием вершин}.
Если граф $\Delta$ можно получить из $\Gamma$ применив слияние вершин несколько раз, то мы будем писать $\Delta\prec \Gamma$.

Заметим, что вышесказанное влечёт следующие два наблюдения:
\begin{itemize}
 \item Если $\Delta$ --- индуцированный подграф конечного графа $\Gamma$ тогда $\Delta\prec \Gamma$.
 \item Если $\Delta\prec \Gamma$, то $\Gamma$-сравнение влечёт $\Delta$-сравнение.
\end{itemize}
Отсюда получаем следующую переформулировку теоремы.

\begin{thm}{Переформулировка}
Для любого конечного связного графа $\Gamma$, 
\begin{enumerate}[(a)]
\item $\Gamma$-сравнение тривиально или
\item $C_4\prec \Gamma$ или
\item $T_3\prec \Gamma$.
\end{enumerate}
\end{thm}

Связный граф будем называть \emph{мультипутём}, если существует целозначная функция $\ell$ определённая на его вершинах такая, что 
$$v\sim w
\quad\Longleftrightarrow\quad
|\ell(v)-\ell(w)|\leqslant 1.$$
Значение $\ell(w)$ будет называться \emph{уровнем} вершины $w$.
Мультипуть полностью определяется целочисленной последовательностью, задающей количество вершин на каждом уровне.
\begin{figure}[ht!]
\centering
\medskip
\includegraphics{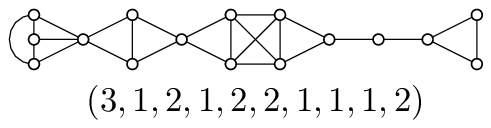}
\medskip
\end{figure}
Пример мультипути и его последовательность показаны на рисунке.

\begin{thm}{Лемма}
Если $\Gamma$ --- конечный связный граф такой, что $C_4\nprec\Gamma$ и $T_3\nprec\Gamma$,
то $\Gamma$ --- мультипуть.
\end{thm}

\parit{Доказательство.}
Обозначим через $|\ -\ |_\Gamma$ метрику кратчайшего пути на~$\Gamma$;
она выдаёт число рёбер в кратчайшем пути между двумя вершинами.
Докажем, что
\[|u-w|_\Gamma\geqslant |u-v|_\Gamma\geqslant|v-w|_\Gamma\geqslant 2
\quad\Longrightarrow\quad |u-w|_\Gamma=|u-v|_\Gamma+|v-w|_\Gamma
\leqno({*})\]
для любых трёх вершин $u$, $v$ и $w$ в $\Gamma$.

Предположим, что $({*})$ не выполняется для $u$, $v$ и $w$.
Рассмотрим подграф $\Delta$ в $\Gamma$, индуцированный тройкой кратчайших путей между $u$, $v$ и $w$.
\begin{figure}[ht!]
\centering
\includegraphics{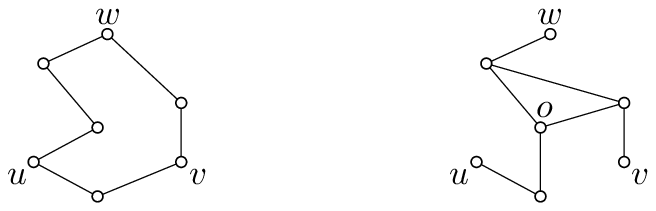}
\end{figure}
Заметим, что $\Delta$ либо является циклом,
либо он содержит три пути из некоторой вершины, скажем $o$, до  $u$, $v$ и $w$ таких, что каждый путь не заходит в оставшиеся вершины тройки.
В этих случаях имеем $C_4\prec\Delta$ или $T_3\prec\Delta$ соответственно.
Применив наблюдение выше, получаем $\Delta\prec\Gamma$ --- $({*})$ доказано.

Обозначим через $d$  диаметр $\Gamma$.
Можно предположить, что $d\geqslant 2$.
Если $d=1$, то $\Gamma$ --- полный граф, в частности, мультипуть.

Выберем вершины $p$ и $q$ такие, что $|p-q|_\Gamma=d$.
Покажем, что $\Gamma$ есть мультипуть со следующей функцией уровня
\[\ell(w)=
\begin{cases}
\hfil |p-w|_\Gamma&\text{если}\quad |p-w|_\Gamma\geqslant 2,
\\
d-|q-w|_\Gamma&\text{если}\quad |q-w|_\Gamma\geqslant 2,
\\
\hfil 1&\hfil\text{в остальных случаях.}
\end{cases}
\]
Из $({*})$ получаем, что $|p-w|_\Gamma+|q-w|_\Gamma=d$, если $|p-w|_\Gamma\geqslant 2$ и $|q-w|_\Gamma\geqslant 2$, а значит $\ell$ хорошо определена.

Если $d\geqslant 4$, то наше утверждение следует из $({*})$.
Остаётся два случая: $d=2$ и $d=3$.

\begin{wrapfigure}{r}{17 mm}
\vskip-6mm
\centering
\includegraphics{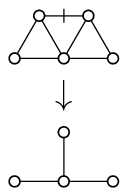}
\vskip-2mm
\end{wrapfigure}

Заметим, что граф-веер на рисунке не может быть индуцированным подграфом $\Gamma$.
Действительно, применяя слияние вершин для концов отмеченного ребра можно получить треногу --- противоречие.

\parit{Случай $d=2$.}
По $({*})$, $\ell^{-1}(0)$ и $\ell^{-1}(2)$ образуют клики.
Заметим, что $\ell(v)=1$ тогда и только тогда, когда $p\sim v$ и $q\sim v$.
Далее заметим, что $\ell^{-1}(1)$ образует клику;
действительно, если $u\nsim v$ для некоторых $u,v\in\ell^{-1}(1)$,
то подграф индуцированный $\{p,q,u,v\}$ есть четыре-цикл --- противоречие.

Остаётся доказать, что $u\sim v$, $v\sim w$, и $u\nsim w$ если $\ell(u)=0$, $\ell(v)=1$, и $\ell(w)=2$.

\begin{wrapfigure}{r}{17 mm}
\vskip-2mm
\centering
\includegraphics{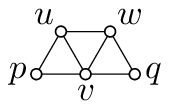}
\vskip4mm
\includegraphics{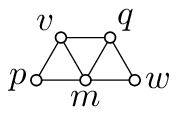}
\end{wrapfigure}

Предположим $u\sim w$.
Заметим, что $v\sim u$ и $v\sim w$, иначе $\Gamma$ содержит индуцированный цикл длины 4 или 5 с вершинами $p$, $u$, $v$, $w$, $q$.
Таким образом, подграф индуцированный $\{p,u,v,w,q\}$ изоморфен вееру --- противоречие.

Предположим $v\nsim w$.
Заметим, что $w\ne q$.
Обозначим через $m$ вершину между $w$ и $p$.
Из вышесказанного, $\ell(m)=1$;
в частности $v\sim m\sim q$.
И снова, подграф индуцированный $\{p,v,m,w,q\}$ изоморфен вееру  --- противоречие.
Точно так же доказывается, что $u\sim v$.

\parit{Случай $d=3$.}
Нам нужно доказать что $u\sim v$ если $\ell(u)=2$ и $\ell(v)=3$;
остальное следует из $({*})$.

\begin{wrapfigure}{r}{20 mm}
\vskip-0mm
\centering
\includegraphics{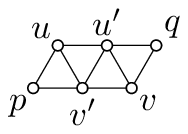}
\vskip-2mm
\end{wrapfigure}

Предположим противное.
Пусть $u'$ ($v'$) лежит посередине $u$ и $q$ (соответственно, $v$ и $p$).
Заметим, что $\ell(u')=3$, $\ell(v')=2$, и $u'\sim v'$ (иначе подграф, индуцированный $\{u,v,u',v'\}$, образует 4-цикл).
То есть, подграф, индуцированный $\{p,q,u,v,u',v'\}$, изображён на рисунке.
Заметим, что он содержит веер --- противоречие.
\qeds

\begin{thm}{Предложение}
Пусть $\Gamma$ --- мультипуть с последовательностью $(k_0,\z\dots, k_m)$.
Предположим $C_4\nprec\Gamma$ и $T_3\nprec\Gamma$.
Тогда
\begin{enumerate}[(a)]
 \item\label{lem:multipath:5} Если $m\geqslant 4$, то $k_2=\dots=k_{m-2}=1$.
 \item\label{lem:multipath:4} Если $m= 3$, то $k_1=1$ или $k_2=1$.
 \item\label{lem:multipath:3} Если $m= 2$, то $k_0=1$, $k_1=1$ или $k_2=1$.
\end{enumerate}

\end{thm}

\parit{Доказательство.} Предполагая противное в каждом случае, получаем

\parit{(\ref{lem:multipath:5})}
если $m\geqslant 4$, то мултипуть $(1,1,2,1,1)$ --- индуцированный подграф в~$\Gamma$,

\parit{(\ref{lem:multipath:4})}
если $m=3$, то мультипуть $(1,2,2,1)$ --- индуцированный подграф в~$\Gamma$,

\parit{(\ref{lem:multipath:3})}
если $m=2$, то мультипуть $(2,2,2)$ --- индуцированный подграф в~$\Gamma$.

В каждом случае мы приходим к противоречию, применив слияние вершин к концам отмеченный рёбер на рисунке.

\begin{figure}[h!]
\centering
\includegraphics{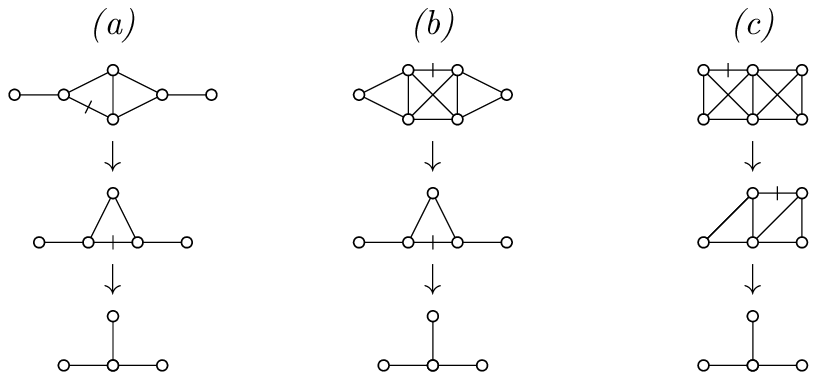}
\end{figure}
\qedsf

Остаётся показать, что для любого мультипути $\Gamma$, описанного в предложении,
$\Gamma$-сравнение тривиально.
Мы проделаем это явно, указав координаты конфигурации сравнения на вещественной прямой.

Каждое ребро в $\Gamma$ имеет вес --- расстояние между его концами в $X$.
Определим расстояние $\|v-w\|_\Gamma$ как минимальный вес пути, соединяющего $v$ и $w$ в $\Gamma$.
Заметим, что 
\[\|v-w\|_\Gamma\geqslant |v-w|_X\]
для любых $v$ и $w$.

Если $m\leqslant 1$, то $\Gamma$ --- полный граф.
В этом случае $\Gamma$-сравнение тривиально.
Остаётся рассмотреть случаи $m\geqslant2$.

Выберем вершину $w$, которая единственна на своём уровне и не очень далека от середины $\Gamma$.
А именно, если $m\geqslant 4$, то возьмём за $w$ вершину на втором уровне; 
согласно предложению, она в нём единственна.
Если $m=3$, то согласно предложению можно предположить, что $k_2=1$;
в этом случае выберем $w$ на втором уровне.
И наконец, если $m=2$, пусть $w$ будет единственной вершиной в своём уровне;
по предложению такая существует.

Для каждой вершины $v_i$, положим
\[\tilde v_i=\pm \|w-v_i\|_\Gamma,\]
здесь стоит знак плюс, если уровень $v_i$ выше уровня $w$, и минус в противном случае.
По неравенству треугольника, полученная конфигурация $\tilde v_1,\dots,\tilde v_n\in\mathbb{R}$ удовлетворяет $\Gamma$-сравнению.
\qeds

\parbf{Замечания.}
Утверждения, перечисленные в предисловии, говорят, что сравнение для правильно выбранного графа (или семейства графов) отвечает за содержательное семейство метрических пространств.
Далее мы приведём ещё пару наблюдений о граф-сравнении общего характера.

Граф-сравнения для всех полных двудольных графов влекут так называемые \emph{чистые неравенства отрицательного типа} \cite[6.1.1]{deza-laurent-ru}.
По критерию Шёнберга, этих неравенств достаточно для существования изометрического вложения в гильбертово пространство \cite[6.2.1]{deza-laurent-ru}.
В частности выполнение сравнений для всех графов влечёт, что метрическое пространство изометрично подмножеству гильбертова пространства.
Последнее утверждение можно доказать на прямую тем же путём, что \cite[Prop. 1.9]{toyoda}.

\begin{wrapfigure}{r}{17 mm}
\vskip-6mm
\centering
\includegraphics{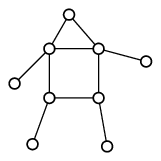}
\vskip-2mm
\end{wrapfigure}

Предыдущие наблюдения выполняются для произвольных метрик.
Для пространств с внутренней метрикой большинство граф-сравнений влекут то, что пространство изометрично подмножеству гильбертова пространства.
Действительно, если $C_4\prec\Gamma$ и $T_3\prec\Gamma$ (как для графа на рисунке),
то любое полное пространство с внутренней метрикой имеет нулевую кривизну в смысле Александрова; 
в частности изометрично выпуклому множеству гильбертова пространства.

\parbf{Благодарности.}
Мы благодарны А. Б. Лычаку за помощь. 

Первый автор частично поддержан Российским фондом фундаментальных исследований, грант 20-01-00070;
второй автор поддержан Национальным научным фондом, грант DMS-2005279
и
Министерством образования и науки Российской Федерации, грант 075-15-2022-289.

{\sloppy
\printbibliography[heading=bibintoc]
\fussy
}

{\bigskip\footnotesize

\noindent Н. Д. Лебедева,
\par\nopagebreak
 \textsc{Санкт-Петербургский государственный университет, Университетская наб., д. 7–9, Санкт-Петербург, 199034, Россия.}
\par
\nopagebreak
 \textsc{Федеральное  государственное  бюджетное  учреждение  науки
Санкт-Петербургское  отделение  Математического  института
им. В. А. Стеклова  Российской  академии  наук, наб. Фонтанки 27,
191023 Санкт-Петербург,
Россия }
  \par\nopagebreak
  \textit{Электронный адрес}: \texttt{lebed@pdmi.ras.ru}

\medskip

\noindent   А. М. Петрунин, 
\par\nopagebreak
 \textsc{Math. Dept. PSU, University Park, PA 16802, USA.}
  \par\nopagebreak
  \textit{Электронный адрес}: \texttt{petrunin@math.psu.edu}
  
}